%%%%%%%%%%%%%%%%%%%%%%%%%%%%%%%%%%%%%%%%%%%%%%%%%%%%%%%%%%%%%%%%%%%%%%%%%%%
%% Gardner, R. J.
%% 
%% On the Busemann-Petty problem concerning central sections of centrally 
%%   symmetric convex bodies
%% 
%% We present a method which shows that in $\Eb$ the Busemann-Petty 
%%   problem, concerning central sections of centrally symmetric convex 
%%   bodies, has a positive answer. Together with other results, this 
%%   settles the problem in each dimension.
%% 
%% publ:  Bull. Amer. Math. Soc. (N.S.) 30(1994) no. 2
%% pp:    222-226
%% type:  Research Announcement        markup: amstex    file size: 21K
%% contact:gardner@baker.math.wwu.edu
%% 
%% copyright: American Math. Society copyright; see end of article
%% 
%% Include files necessary for this article: bull-ppt.tex
%% 
%%%%%%%%%%%%%%%%%%%%%%%%%%%%%%%%%%%%%%%%%%%%%%%%%%%%%%%%%%%%%%%%%%%%%%%%%%%
\input amstex
\documentstyle{amsppt}
\input bull-ppt
\keyedby{bull493/lic}

%\magnification=\magstep1

\define \En{\Bbb{E}^n}

\define \Eb{\Bbb{E}^3}
\define \Ec{\Bbb{E}^4}

\define \Rn{\Bbb{R}^n}

\predefine\preres{\restriction}
\redefine\restriction{|}

\predefine\preconv{\conv}
\redefine \conv{\operatorname{conv}}

\predefine\accc{\c}
\redefine   \c{\operatorname{c}}

%\redefine   \log{\operatorname{log}}
%\redefine   \exp{\operatorname{exp}}

%\define   \supp{\operatorname{supp}}

\topmatter
\cvol{30}
\cvolyear{1994}
\cmonth{April}
\cyear{1994}
\cvolno{2}
\cpgs{222-226}
\ratitle
\title 
On the Busemann-Petty problem\\ concerning central 
sections\\ of  
centrally symmetric convex bodies
\endtitle
\author R. J. Gardner 
\endauthor
\shorttitle %\rightheadtext\nofrills{\sevenrm 
{On the Busemann-Petty Problem}
%\leftheadtext\nofrills{\sevenrm R. J. GARDNER}
\address 
Department of Mathematics, Western Washington University, 
Bellingham,
Washington 98225-9063 \endaddress
%\endaddress
%\email
\ml gardner{\@}baker.math.wwu.edu\endml
%\endemail
\thanks Supported in part by NSF Grant DMS 9201508\endthanks
 %\endsubjclass
\date May 11, 1993\enddate
\abstract We present a method which shows that in $\Eb$ the 
Busemann-Petty  
problem, concerning central sections of centrally 
symmetric convex  
bodies, has a positive answer.  Together with other 
results, this  
settles the problem in each dimension.\endabstract 
\endtopmatter

\document

%\NoBlackBoxes

In \cite{BP}, Busemann and Petty asked the following 
question,
which resulted from reformulating a problem in Minkowskian 
geometry.   
Suppose $K_1$ and $K_2$ are convex bodies in 
$n$-dimensional  
Euclidean space $\En$ which are centered (centrally 
symmetric with  
center at the origin) and such that
$${\lambda}_{n-1}(K_1\cap u^{\perp})\le 
{\lambda}_{n-1}(K_2\cap  
u^{\perp})\,,$$
for all $u$ in the unit sphere $S^{n-1}$.  Then is it true 
that
$$\lambda_n(K_1)\le \lambda_n(K_2)\,?$$
(Here $u^{\perp}$ denotes the hyperplane through the 
origin  
orthogonal to $u$, and $\lambda_k$ denotes $k$-dimensional 
Lebesgue  
measure, which we identify throughout with $k$-dimensional 
Hausdorff  
measure.)

The question, now generally known as the Busemann-Petty 
problem, has  
often appeared in the literature.  More than thirty years 
ago, Busemann  
gave the problem wide exposure in \cite{B2} and Klee posed  
it again in \cite{K}.  The problem attracted the attention 
of  
those working in the local theory of Banach spaces; see, 
for example,  
the paper \cite{MP, p. 99} of Milman and Pajor.  It 
surfaces  
again in Berger's article \cite{Be, p. 663}, and it is 
also  
stated in the books of Burago and Zalgaller \cite{BZ, p. 
154};
Croft, Falconer, and Guy \cite{CFG, Problem A9, p. 22}; 
and  
Schneider \cite{S, p. 423}.

The problem has an interesting history.  Using a clever 
probabilistic  
argument, Larman and Rogers \cite{LR} proved that the 
answer,  
surprisingly, is negative in $\En$ for $n\ge 12$.  Later, 
Ball  
\cite{B} applied his work on maximal sections of a cube to 
obtain  
a negative answer for $n\ge 10$, where $K_1$ is a centered 
cube and  
$K_2$ a centered ball of suitable radius.  Giannapoulos 
\cite{Gi}  
improved this negative result to $n\ge 7$ by using an 
appropriate  
cylinder for $K_1$ instead of a cube.  Independently, 
Bourgain  
\cite{Bo} showed that the same result can be achieved by 
taking  
$K_1$ to be a suitable arbitrarily small perturbation of a 
centered  
ball; Bourgain also proved that his method will not work 
in $\Eb$.  A  
further improvement was made by Papadimitrakis \cite{P} 
and the  
author \cite{G1}, independently, by demonstrating that the  
answer is negative for $n\ge 5$, when $K_1$ is a centered 
cylinder.   
The most recent negative answer was obtained, for $n\ge 4$ 
and $K_1$  
a centered cube, by Zhang \cite{Z1}, \cite{Z2}.   
For $4\le n\le 6$, the existence of a suitable $K_2$ can 
be proved,  
though it seems likely that for these values of $n$, $K_2$ 
cannot be  
a ball, regardless of the choice of $K_1$.  Other papers 
on the  
problem related to those mentioned above include \cite{GR} 
and 
\cite{T}.

We shall outline a solution for $n=3$, thereby settling 
the problem  
in each dimension.  Against the background of the results 
above, the  
positive answer for $n=3$ is unexpected.  It is also 
especially  
interesting from the point of view of geometric 
tomography, in which  
one attempts to obtain information about a geometric 
object from data  
concerning its sections or projections.  Geometric 
tomography has  
connections with functional analysis and possible 
applications to  
robot vision and stereology (see, for example, \cite{BL,  
ES, GV, MP, W}).

A few positive results are already known.  The case $n=2$ 
is trivial.   
Busemann and Petty themselves noted in \cite{BP} that the  
Busemann intersection inequality (see \cite{B1, (4), p. 
2})  
may be applied to obtain a positive answer when $K_1$ is a 
centered  
ellipsoid.  Lutwak \cite{L} obtained an important 
generalization  
of this fact by showing that the same is true whenever 
$K_1$ is a  
member of a certain class of bodies which he called 
intersection  
bodies.

We shall explain Lutwak's result in some detail, since it 
is an  
essential ingredient in our method.  A set $L$ in $\En$ is 
{\it  
star shaped} at the origin if it contains the origin and 
if every line  
through the origin meets $L$ in a (possibly degenerate) 
line segment.   
By a {\it star body} we mean a compact set $L$ which is 
star shaped  
at the origin and whose radial function
$$\rho_L(u)=\text {max}\{c\ge 0:cu\in L\}\,,$$
for $u\in S^{n-1}$, is continuous on $S^{n-1}$.  The star 
body $L$ is  
called the {\it intersection body} of another star body 
$M$ if
$$\rho_L(u)={\lambda}_{n-1}(M\cap u^{\perp})\,,$$
for all $u\in S^{n-1}$.  We write $L=IM$; it is clear that 
$L$ must be  
centered, and it is known (see \cite{L}) that there is a 
unique  
centered star body $M'$ for which $L=IM'$.

A useful alternative viewpoint is provided by the 
spherical Radon  
transform.  
Suppose $g$ is a Borel function on $S^{n-1}$ and $f$ is 
defined by
$$f(u)=\int_{S^{n-1}\cap 
u^{\perp}}g(v)\,d\lambda_{n-2}(v)\,,$$
for all $u\in S^{n-1}$; that is, $f(u)$ is the integral of 
$g$ over  
the great sphere in $S^{n-1}$ orthogonal to $u$.  Then we 
write
$$f=Rg$$
and say that $f$ is the {\it spherical Radon transform} of 
$g$.   
Using the polar coordinate formula for volume, we see that 
a star  
body $L$ is the intersection body of some star body $M$ if 
and only  
if $\rho_L=Rg$ for some nonnegative continuous function 
$g$; just  
take $g=\rho_M^{n-1}/(n-1)$.

Suppose $L_1$ is the intersection body of some star body 
and $L_2$ is  
an arbitrary star body, such that   
$${\lambda}_{n-1}(L_1\cap u^{\perp})\le 
{\lambda}_{n-1}(L_2\cap  
u^{\perp})\,,$$
for all $u\in S^{n-1}$.  Then Lutwak's theorem (see 
\cite{L,  
Theorem 10.1}) says that $\lambda_n(L_1)\le \lambda_n(L_2)$.

It is worth noting that Lutwak has offered as an 
alternative and  
different definition of intersection body, a star body $L$ 
such that  
$\rho_{L}=R\mu$, where $\mu$ is an even finite Borel 
measure in  
$S^{n-1}$.  A consequence of a result in \cite{GLW} is 
that  
Lutwak's theorem still holds when $L_1$ is an intersection 
body in  
this wider sense of the term.

The class of intersection bodies is in a sense dual to the 
better-known class of projection bodies.  The latter, 
which are just the  
centered zonoids, have been intensively studied and have 
many  
applications; see, for example, the articles of Bourgain 
and  
Lindenstrauss \cite{BL}, Goodey and Weil \cite{GW}, and  
Schneider and Weil \cite{SW}, or Schneider's book \cite{S,  
Section 3.5}.  In fact, the Busemann-Petty problem has a 
dual form in  
which sections are replaced by projections.  This dual 
problem was  
solved, by Petty and Schneider independently, shortly 
after it was  
posed; the answer is negative for all $n>2$.  Lutwak's 
theorem  
concerning sections of star bodies is also dual to a 
corresponding  
one for projections of convex bodies, obtained by Petty 
and Schneider  
using tools from the Brunn-Minkowski theory (see \cite{S, 
p.  
422}).  For sections, the extension from convex bodies to 
star bodies  
is not only natural but crucial.  For example, it can be 
seen by  
direct calculation that a centered cylinder in $\Eb$ is 
the  
intersection body of a nonconvex centered star body; see 
\cite{G1, 
Remark 5.2(ii)}.

The other known positive results are as follows.  Hadwiger 
\cite{H} 
and Giertz \cite{Gie} independently showed that the 
question  
has an affirmative answer when $K_1$ and $K_2$ are coaxial 
centered  
convex bodies of revolution in $\Eb$.  In \cite{G1,  
Theorem 5.1}, it is proved that a centered convex body of 
revolution  
$K$ whose radial function $\rho_K$ belongs to  
$C^{\infty}_e(S^{n-1})$, the class of infinitely 
differentiable even  
functions on the unit sphere, is the intersection body of 
some star  
body when $n=3$ or $4$.  (This result is re-proved in 
\cite{Z2} 
by a different method.)  Using Lutwak's theorem, it is 
easy  
to see that this implies that the Busemann-Petty problem 
has a  
positive answer in $\Eb$ or $\Ec$ whenever $K_1$ is a 
centered convex  
body of revolution.  It has also been shown by Meyer 
\cite{M}  
that the answer is positive in $\En$ provided that $K_1$ 
is a  
centered cross-polytope (the $n$-dimensional version of 
the  
octahedron).

We now sketch a proof that the Busemann-Petty problem has 
a positive  
answer in $\Eb$.  The details will appear in \cite{G2}.

\proclaim{Theorem}
The Busemann-Petty problem has an affirmative answer in 
$\Eb$.
\endproclaim

\demo{Sketch of the proof}
By approximating, we can assume that $\rho_{K_1}\in  
C^{\infty}_e(S^{2})$ and that $K_1$ has everywhere 
positive Gaussian  
curvature.  By Lutwak's theorem, it suffices to show that 
an  
arbitrary centered convex body $K$ in $\Eb$ with these 
additional  
properties is the intersection body of some star body.  
(We make no  
attempt in this note to find the least restrictive 
additional  
conditions to impose on $K$.)  This will be proved if 
there is a  
nonnegative function $g\in C(S^{2})$ such that 
$\rho_K=Rg$.  It is  
known that since\ $\rho_{K}\in C^{\infty}_e(S^{2})$, a 
$g\in  
C^{\infty}_e(S^{2})$ exists and is unique.  Let $u_0\in 
S^{2}$.  An  
inversion formula of Funk \cite{F} gives
$$g(u_0)=\text {lim}_{t\rightarrow  
1^{-}}\frac{1}{2\pi}\frac{d}{dt}\int_0^t\frac{xA_K(\text  
{sin}^{-1}x)}{\sqrt{t^2-x^2}}\,dx\,,$$
where $A_K(\phi)$ denotes the average of $\rho_K$ on the 
circle of  
latitude with angle $\phi$ from the north pole $u_0$.  
After some  
manipulation, one can obtain
$$2\pi  
g(u_0)=\rho_K(u_0)+
\frac{1}{2\pi}\int_0^{2\pi}\int_0^{\frac{\pi}{2}}
\frac{\partial  
\rho_K(\theta,\phi)}{\partial\phi}\sec\phi\,d\phi\,d\theta%
\,,\tag1$$
where $(\theta,\phi)$ denotes the usual angles of 
spherical polar  
coordinates.

%>
From $K$, construct a body $\bar{K}$, called a Schwarz 
symmetral of  
$K$, as follows.  Each horizontal section of $\bar{K}$ is 
a disk  
whose center lies on the $z$-axis and whose area equals 
that of the  
horizontal section of $K$ of the same height.  From the  
Brunn-Minkowski theorem, it follows that $K$ is a convex 
body of  
revolution, and our assumptions about $K$ allow us to 
conclude that  
$\rho_{\bar K}\in C^{1}_e(S^{2})$.  One can then show that  
$\rho_{\bar{K}}=R\bar{g}$, for some $\bar{g}\in C(S^2)$, 
and that  
equation (1) holds when $g$ and $\rho_K$ are replaced by 
$\bar{g}$  
and $\rho_{\bar{K}}$, respectively.  Moreover, the 
argument of  
\cite{G1, Theorem 5.1} proves that $\bar{g}$ is  
nonnegative.

The final step of the proof involves applying a 
cylindrical  
transformation to equation (1).  Once this is done, it can 
be seen  
that $g(u_0)=\bar{g}(u_0)$, and therefore $g(u_0)\ge 0$, 
as required.   
\qed\enddemo

\enddocument